\documentclass{amsart}
\usepackage{amssymb,euscript,amsmath, mathrsfs}
\usepackage[dvips]{graphicx}
\usepackage[dvips]{color}

\newcounter{ENUM}
\newcommand{\itm}{\item}
\newenvironment{ilist}{\renewcommand{\theENUM}{\roman{ENUM}}\renewcommand{\itm}{\addtocounter{ENUM}{1}\item[(\theENUM)]}\begin{itemize}\setcounter{ENUM}{0}}{\end{itemize}}

\newenvironment{alist}[1][0]{\renewcommand{\theENUM}{\alph{ENUM}}\renewcommand{\itm}{\addtocounter{ENUM}{1}\item[\theENUM)]}\begin{itemize}\setcounter{ENUM}{#1}}{\end{itemize}}

\newcommand{\margh}[1]{}

\def\risom{\overset{\sim}{\rightarrow}}

\input xy
\xyoption{all}
\CompileMatrices

\def\Z{{\mathbb Z}}

\def\P{{\mathbb P}}
\def\Q{{\mathbb Q}}

\def\H{{\mathbb H}}

\def\F{{\mathscr F}}
\def\sG{{\mathscr G}}
\def\cH{{\mathscr H}}

\def\L{{\mathscr L}}
\def\M{{\mathscr M}}
\def\E{{\mathscr E}}
\def\O{{\mathscr O}}

\def\K{{\mathscr K}}
\def\sA{{\mathscr A}}

\def\cM{{\mathcal M}}
\def\T{{\mathcal T}}
\def\U{{\mathcal U}}
\def\cHom{\mathcal{H}om}
\def\cEnd{\mathcal{E}nd}

\def\m{{\mathfrak m}}
\def\n{{\mathfrak n}}

\def\vp{\varphi}

\def\length{\operatorname{length}}
\def\ch{\operatorname{char}}

\def\Ext{\operatorname{Ext}}

\def\Hom{\operatorname{Hom}}

\def\PGL{\operatorname{PGL}}

\def\Spec{\operatorname{Spec}}
\def\Pic{\operatorname{Pic}}

\def\Quot{\operatorname{Quot}}

\def\Gr{\operatorname{Gr}}
\def\Def{\operatorname{Def}}

\def\coker{\operatorname{coker}}
\def\rk{\operatorname{rk}}
\def\im{\operatorname{im}}

\def\red{\operatorname{red}}

\def\oH{\operatorname{H}}

\numberwithin{equation}{section}
\newtheorem{thm}{Theorem}[section]
\newtheorem{prop}[thm]{Proposition}
\newtheorem{lem}[thm]{Lemma}
\newtheorem{cor}[thm]{Corollary}

\theoremstyle{definition}
\newtheorem{defn}[thm]{Definition}
\newtheorem{ques}[thm]{Question}
\newtheorem{ex}[thm]{Example}
\newtheorem{sit}[thm]{Situation}
\theoremstyle{remark}
\newtheorem{notn}[thm]{Notation}
\newtheorem{rem}[thm]{Remark}

\begin{document}
\title{The Generalized {V}erschiebung Map for Curves of Genus 2}
\author{Brian Osserman}
\begin{abstract} Let $C$ be a smooth curve, and $M_r(C)$ the coarse moduli space of vector bundles of rank $r$ and trivial determinant on $C$. We examine the generalized Verschiebung map $V_r:M_r(C^{(p)}) \dashrightarrow M_r(C)$ induced by pulling back under Frobenius. Our main result is a computation of the degree of $V_2$ for a general $C$ of genus $2$, in characteristic $p>2$. We also give several general background results on the Verschiebung in an appendix.
\end{abstract}
\thanks{This paper was partially supported by fellowships from the National Science Foundation and Japan Society for the Promotion of Sciences.}
\maketitle

\section{Introduction}

In this paper we address the degree of the Verschiebung rational map 
$V_{2}$ induced by pullback under Frobenius on the moduli space $M_2$ of 
rank $2$
vector bundles with trivial determinant on a smooth proper curve, in the
case of genus $2$, along with some related questions. 
We work throughout over an algebraically closed field $k$ of characteristic $p>2$, except where specified.

Aside from the importance of the Verschiebung in the case of Jacobians, and the consequent desire to understand its generalization to higher rank,
motivation for understanding the 
geometry of the Verschiebung map was provided by the close relationship
between the Verschiebung map and characteristic-$p$ representations of the fundamental group of $C$, when the base field $k$ for our curve $C$ is finite (see the introduction to \cite{l-p}). In particular, A. J. de Jong 
showed that curves in the moduli
space of vector bundles which are fixed under some iterate of the
Verschiebung will correspond to characteristic-$p$ representations for which the
geometric fundamental group has infinite image, which he conjectures in
\cite{dj2} cannot happen for characteristic-$\ell$ representations. He further shows that such curves would have to pass through the undefined locus of the Verschiebung.
Another motivation comes from the fact that invariants such as the degree 
of Verschiebung nearly always seem to be given by polynomials in $p$, with 
no obvious explanation for why this should be the case. One might hope that examination of enough different cases of this phenomenon would ultimately provide insight into the general situation.

In order to state our result, we fix some terminology:

\begin{defn} Let $C$ be a smooth, proper curve. A semistable vector bundle $\F$ on $C^{(p)}$ is said to be {\bf Frobenius-unstable} if $F^* \F$ is unstable. Now suppose $C$ has genus $2$. Given a Frobenius-unstable vector bundle $\F$ of rank $2$ and trivial determinant, we say that $\F$ is a {\bf reduced} point of the Frobenius-unstable locus if the first-order determinant-preserving infinitesmal deformations of $\F$ inject into the first-order infinitesmal deformations of $F^* \F$.
\end{defn} 

We will show:

\begin{thm}\label{deg-main} Let $C$ be a smooth, proper genus $2$ curve
over an algebraically closed field $k$ of characteristic $p>2$, and suppose
that the Frobenius-unstable locus for vector bundles of rank $2$ and trivial
determinant is composed of $\delta$ reduced points. Then:
\begin{ilist}
\itm Each undefined point of $V_{2}$ may be resolved by a single
blowup, and $V_{2}$ has degree $p^3-\delta$;
\itm The exceptional divisor associated to such an undefined
point maps bijectively to $\P\Ext(\L, \L^{-1}) \subset M_{2}(C)$, where $\L$ is a theta
characteristic on $C$, and specifically is the destabilizing line bundle for
$F^* \F$, where $\F$ is the Frobenius-unstable vector bundle associated to
the undefined point.
\end{ilist}
\end{thm}

It follows from results of Mochizuki (see \cite{mo3}, \cite{os6}) or of the
author (see \cite{os10}) that the hypothesis of Theorem \ref{deg-main} 
holds with $\delta=\frac{2}{3}(p^3-p)$, so we can conclude the following.

\begin{thm} With the notation of Theorem \ref{deg-main}, if $C$ is general, 
we have that $\deg V_2 = \frac{p^3+2p}{3}$. We further know that the 
undefined points of $V_2$ may each be resolved by a single blowup, with 
exceptional divisor mapping bijectively to $\P(\Ext(\L,\L^{-1})$.
\end{thm}

We begin in Section \ref{s-deg-prelim} with some straightforward and general
results on degrees of rational maps of projective spaces, and their
application to $V_{2}$. We next carry out some calculations involving
deformation theory of bundles with connection in Section \ref{s-deg-hyper}.
Finally, in Section \ref{s-deg-geometric} we establish, following the
argument of Langton's properness theorem, the necessary relationship between
the locus of $p$-curvature $0$ 
connections on the unstable bundles $\E$ of \cite[Prop. 2.6]{os11}, 
the undefined locus of $V_2$, 
and the image of the exceptional divisor if one blows up to resolve $V_2$. 
Appendix \ref{s-deg-background} is a compilation of necessary
technical background results on $V_{n}$, presented in more generality and including arguments suggested by A. J. de Jong and Christian Pauly.
Finally, Appendix
\ref{s-deg-digress} develops some slightly non-standard commutative algebra for non-reduced rings which arises in our vector bundle manipulations.

The existing literature on such geometric questions on the Verschiebung 
is considerably scarcer than on Frobenius-unstable vector bundles. The only
other such results in this area were developed recently by Laszlo
and Pauly, who gave explicit polynomials defining the Verschiebung
in the particular cases of genus 2, rank 2, and characteristics 2 and 3, in
\cite{l-p} and \cite{l-p2}. Lange and Pauly also obtain our formula for the
degree of $V_{2}$ via a different approach in \cite{l-p3}, although
their techniques thus far gives only that it is an upper bound in the case 
of ordinary curves, and not that it is an equality. 

The contents of this paper form a portion of the author's 2004 PhD thesis 
at MIT, under the direction of Johan de Jong.

\section*{Acknowledgements}

I would like to thank Johan de Jong for his tireless and invaluable
guidance. I would also like to thank David Sheppard, David Helm, Brian
Conrad, Christian Pauly, and Yves Laszlo for their helpful conversations.

\section{Preliminaries on Degrees}\label{s-deg-prelim}

In this section we make some basic observations about degrees of rational
maps from projective space to itself, and then apply these to the specific case of the Verschiebung. We remark that unless otherwise
specified, in this section the term `point' shall always refer to a closed 
point. We suppose we are in the following situation.

\begin{sit}\label{deg-proj-sit} We are given a rational map 
$f:\P^n \rightarrow \P^n$ which is dominant and defined at all but a finite 
set of points. We suppose that we are given homogeneous coordinates $X_i$ 
on $\P^n$, and $f$ is represented by $n+1$ homogeneous polynomials $F_i$ of 
degree $d$ in the $X_i$.
\end{sit}

Although special cases of the following proposition are certainly extremely well-known, the general statement, and in particular the possibility of inequality, appears less widely known.

\begin{prop}\label{deg-formula} In the above situation, we have the
inequality:
$$\deg f \leq d^n - \delta,$$ where $\delta$ is the total length of the
`undefined locus' subscheme $E_f$ of $\P^n$ cut out by the $F_i$.
Moreover, the following are equivalent:
\begin{alist}
\itm The above inequality is an equality;
\itm $E_f$ is a locally complete intersection; 
\itm $E_f$ is Gorenstein. 
\end{alist}

In particular, we get equality when the length of the points of $E_f$ 
are all $1$ or $2$.
\end{prop}

\begin{proof}
Choose $P$ in $\P^n$, and write $H_1, \dots, H_n$ for $n$ hyperplanes 
cutting out $P$. Denote by $E_P$ the (scheme-theoretic) intersection
of the $f^*(H_i)$. Now, for a given $X_{i}$ on the image
space, we observe that on $\P^n_{f^* X_{i}}=\P^n_{F_{i}}$, we have
$f^{-1} (P) \cong E_P$. As ${i}$ varies, the $\P^n_{F_{i}}$ will cover 
everything except $E_f$. Since the $f^*(H_i)$ are given by polynomials homogeneous and linear in the $F_i$, there is a closed immersion of $E_f$
into $E_P$. Hence, we can write $E_P$ as a set as $E_f \cup f^{-1}(P)$,
with a closed immersion of the latter into the former. For $P$ general, $f^{-1}(P)$ and hence $E_P$ is $0$-dimensional, with $\length f^{-1}(P)=\deg f$. Since $\delta:=\length E_f$, and $E_P$ is a complete intersection of $n$ hypersurfaces of degree $d$, Bezout's theorem yields the desired inequality.

Now, $E_f$ being a locally complete intersection implies that it is 
Gorenstein, and conversely, since $E_f$ is cut out by $n+1$ hypersurfaces, if it is Gorenstein it must be a local complete intersection; see \cite[Cor 21.19]{ei1} and \cite[p. 542]{ei1}. Next, to see that a) implies b), we note if we have equality, $E_f \rightarrow E_P$ is an isomorphism onto its image, so since $E_P$ is a complete intersection, $E_f$ is a locally complete intersection.

Finally, we need to show that if $E_f$ 
is a locally complete intersection, the ideal generated by the 
$f^*(H_i)$ corresponding to a general point $P$ on the image $\P^n$ is 
locally equal to the ideal of the $F_i$ at each of the finitely many $Q$ of
$E_f$; we can therefore check the statement one $Q$ at a time. We fix some $Q \in E_f$, and assume we have 
fixed a choice of dehomogenization, so that the $F_i$ are actually functions 
on a neighborhood of $Q$. We also note that it will be enough to prove the
statement for any particular choice of $H_i$ cutting out the given $P$. We
first observe that since $E_f$ is locally a complete intersection,
Nakayama's lemma implies that the defining ideal $I_{E_f}$ of $E_f$ at $Q$ may be generated by $n$ of the $n+1$ generators $F_i$. Hence, we have 
$F_j = \sum_{i \neq j} a _i F_i$ for some integer $j$ and 
$a_i \in \O_{\P^n, Q}$; reindexing, we can assume $j=n+1$. Now, a general point $P$ in the image $\P^n$ may be 
cut out by $H_i$ of the form $X_i - \lambda _i X_{n+1}$ for $i\leq n$. One then checks directly that for a general choice of the $\lambda_i$, $F_{n+1}$ is in the ideal generated by the $f^*(H_i)$. It follows that the other $F_i$ will also be in the ideal generated by the $f^*(H_i)$, giving $E_f \cong E_P$ at $Q$, as desired.
\end{proof}

For posterity, we observe that the inequality of the preceding
proposition need not be an equality:

\begin{ex}Consider the map from $\P^2$ to itself given by 
$(X,Y,Z)\mapsto (X^3, Y^3, XYZ)$. This is undefined only at $(0,0,1)$,
where the subscheme $E_f$ has length $5$, and is visibly not a local complete intersection. We can check that the map is
dominant of degree $3$, so as required, we have strict inequality in the preceding proposition.
\end{ex}

We now apply these results to examine what we can say about the degree of
the Verschiebung map induced by pullback under Frobenius on the moduli space
of vector bundles of rank 2 with trivial determinant on a curve $C$ of genus
$2$. We begin by reviewing the basic facts about the moduli space and
Verschiebung, without any hypotheses on the genus. 

We will denote by $M_2$ the moduli space of semistable vector bundles
of rank $2$ and trivial determinant on $C$, which will be the only space
we will consider here. We will refer to the map induced by pullback under 
Frobenius on moduli spaces of vector bundles as the {\bf Verschiebung},
since in the case of line bundles this is precisely what it is. 
We denote by $V_{2}$ the particular Verschiebung map on $M_{2}$.
$V_{2}$ is a dominant rational map, with finitely many undefined points
composed of the Frobenius-unstable vector bundles. See Appendix
\ref{s-deg-background} for the precise technical definitions and proofs of 
these statements, as well as \cite[Thm. 3.2]{j-x} for the finiteness result.

Now, in our situation of rank $2$ bundles on a curve of genus $2$, 
it is a theorem of Narasimhan and Ramanan that $M_{2} \cong \P^3$; see 
\cite[Thm. 2, \S 7]{n-r}, and note that despite the Riemann surface 
language, the argument goes through unmodified in arbitrary odd 
characteristic. We also know that the Verschiebung map is given by polynomials of degree $p$; see \cite[Prop. A.2]{l-p2}, and note that $\O(1)$ on our
$\P^3$ pulls back to the inverse of the determinant bundle on the moduli 
stack, by Lemma \ref{det-o1}. Alternatively, one can argue in our case by 
pulling back to the Jacobians of $C$ and $C^{(p)}$, which map to Kummer 
surfaces inside $M_{2} \cong \P^3$.

Putting this together with Proposition \ref{deg-formula} yields the following.

\begin{cor}\label{deg-versch} For curves of genus 2, the degree of $V_{2}$
is bounded above by $p^3$; or 
more sharply, $p^3-\delta$, where $\delta$ is the number of points at
which $V_{2}$ it is undefined. If the undefined points are reduced, this
upper bound is an equality.
\end{cor}

On the other hand, we also have:

\begin{lem}\label{deg-lowerbd} For curves of genus $2$, the degree of
$V_{2}$ is bounded below by $p^2$. 
\end{lem}

\begin{proof}For any $\L$ of degree $0$ and not $2$-torsion, consider the bundle $\L \oplus \L^{-1}$. There are $p^2$ line bundles of degree $0$ mapping to $\L$
under $V_{2}$ on the Jacobian, differing from each other by a $p$-torsion
line bundle, and since $\L$ is not $2$-torsion, this gives $p^2$ bundles mapping to $\L \oplus \L^{-1}$. Thus, on an open part of the Kummer surface, each point has at least $p^2$
preimages, and noting that the points having positive-dimensional preimage must have codimension at least $1$ in the Kummer surface, we conclude the desired lower bound on the degree of $V_{2}$.
\end{proof}

\section{Deformation Theory Calculations}\label{s-deg-hyper}

We begin by reviewing some fundamental facts about representing groups of
deformations with cohomology and hypercohomology. Throughout the remainder
of this paper, excepting the appendices, we fix the notation:

\begin{notn}\label{deg-deform} A `deformation' refers to a first order 
infinitesmal deformation, and $\epsilon$ is a square-zero element.
\end{notn}

Let $C$ be any smooth curve, with vector bundle $\E$ and connection $\nabla$. The following propositions are well-known, and may be checked explicitly in terms of Cech cocycles on the $U_i$.

\begin{prop}The space of deformations of $\E$ is isomorphic to
$H^1(C, \cEnd(\E))$, and the space of deformations of $\E$ preserving 
the determinant of $\E$ is isomorphic to $H^1(C, \cEnd^0(\E))$.
\end{prop}

\begin{prop}The space of deformations of
$\nabla$ over $\E$ (respectively, fixing the determinant of $\nabla$)
is isomorphic to $H^0(C, \cEnd(\E)\otimes \Omega^1_C)$ (respectively, 
$H^0(C, \cEnd^0(\E)\otimes \Omega^1_C)$). 
\end{prop}

In order to analyze transport of $\nabla$ along automorphisms of $\E$, we introduce the sheaf map which is described by the connection $\nabla^{\cEnd}$ induced on $\cEnd(\E)$ by $\nabla$; we denote it 
$d_\nabla : \cEnd(\E) \rightarrow \cEnd(\E) \otimes \Omega^1_C$, and it is
given by $\phi \mapsto \nabla \circ \phi - \phi
\circ \nabla$. We then have the following.

\begin{prop}The space of transport equivalence classes of deformations of
$\nabla$ over $\E$ is isomorphic to $H^0(C, \cEnd(\E)\otimes \Omega^1_C)/
d_\nabla (H^0(C, \cEnd(\E)))$. 
Furthermore, $d_\nabla |_{\cEnd^0(\E)}$ takes values in $\cEnd^0(\E)$, and 
if $\ch k$ is prime to the rank of $\E$, the space of transport 
equivalence classes of deformations of $\nabla$ over $\E$ with fixed
determinant is isomorphic to $H^0(C, \cEnd^0(\E)\otimes \Omega^1_C)/
d_\nabla (H^0(C, \cEnd^0(\E)))$. 
\end{prop}

\begin{proof}Beyond the usual cocycle-level verifications, the only 
observation is that since $\ch k$ is prime to the rank, and scalar 
automorphisms fix connections under transport, we have 
$d_\nabla (H^0(C, \cEnd^0(\E))) = d_\nabla(H^0(C, \cEnd(\E)))$.
\end{proof}

However, parametrizing pairs of deformations of $\E$ together with
deformations of $\nabla$ over $\E$ cannot be described naturally with sheaf
cohomology, but rather requires sheaf hypercohomology. Since our complex will have only two terms, hypercohomology in our situation is particularly simple. Given $f: \F_1 \rightarrow \F_2$ an element of $\H^1$ is given by a pair of a Cech 1-cocyle of $\F_1$ and a Cech 0-cochain of $\F_2$ which agree under $f$ and the Cech coboundary map.

\begin{notn} The complex $\cEnd^0(\E) \overset{d_{\nabla}}{\rightarrow} 
\cEnd^0(\E) \otimes \Omega^1_C$ will also be denoted by $\K^{\bullet}$. 
\end{notn}

Given the explicit description of hypercohomology, the following proposition may also be checked directly in terms of our trivializations.

\begin{prop}The space of transport equivalence classes of deformations of 
$\E$ together with $\nabla$
(respectively, deformations fixing both determinants)
is isomorphic to $\H^1(C, \cEnd(\E) \overset{d_{\nabla}}{\rightarrow} 
\cEnd(\E)\otimes \Omega^1_C)$
(respectively, $\H^1(C, \cEnd^0(\E) \overset{d_{\nabla}}{\rightarrow} 
\cEnd^0(\E)\otimes
\Omega^1_C)$).
\end{prop}

Now, every deformation $\nabla'$ of $\nabla$ on $\E$ gives rise to the obvious pair of deformations of $(\E, \nabla)$: namely,
$(\E, \nabla')$. Similarly, a pair $(\E', \nabla')$ yields 
a deformation $\E'$ of $\E$. This gives us induced maps, but we use the standard spectral sequence for the cohomology of a double
complex to show that we obtain a short exact sequence, and obtain some additional information. 

\begin{prop}\label{deg-les} Suppose $\nabla$ has $p$-curvature $0$, and the 
corresponding
$\F$ which pulls back under Frobenius to $\E$ is stable. Then $d_{\nabla}$
is injective on global sections, so we can consider 
$H^0(C, \cEnd^0(\E))$ to be a subgroup of 
$H^0(C, \cEnd^0(\E) \otimes \Omega^1_C)$, and we left get a left exact 
sequence
$$0 \rightarrow H^0 (C, \cEnd^0(\E)\otimes \Omega^1_C)/H^0(C, \cEnd^0(\E)) 
\rightarrow \H^1(C, \K^{\bullet}) \rightarrow H^1(C, \cEnd^0(\E))$$
where the image on the right is the image
under Frobenius pullback of $H^1(C^{(p)}, \cEnd^0(\F))$.
\end{prop}

\begin{proof}Starting with the Cech double complex
associated to our complex,
and taking differentials in the vertical direction, we see that 
the $E_1$ term is:

$$\xymatrix{
{\vdots} & {\vdots} & {} & {}\\
{0} & {0} & {} & {}\\
{H^1(C, \cEnd^0(\E))} \ar[r]^-{d_2^1} & {H^1(C, \cEnd^0(\E) \otimes \Omega^1_C)} 
& {0} & {\dots}\\
{H^0(C, \cEnd^0(\E))} \ar[r]^-{d_2^0} & {H^0(C, \cEnd^0(\E) \otimes \Omega^1_C)}
& {0} & {\dots} }$$
where $d_2^i$ are the maps induced on Cech $i$-cocycles by $d_\nabla$. We 
see immediately that the spectral sequence stabilizes at
$E_2$, and yields a short exact sequence:
$$0 \rightarrow \coker d_2^0 \rightarrow 
\H^1(C, \cEnd^0(\E) \rightarrow \cEnd^0(\E)\otimes \Omega^1_C) 
\rightarrow \ker d_2^1 \rightarrow 0.$$

An element in the kernel of $d_2^0$ is a trace zero
endomorphism of $\E$ which commutes with $\nabla$. But this is precisely
the condition for it to come from a trace zero endomorpism of $\F$ (see,
for instance, \cite[Thm 5.1]{ka1}),
which must be $0$, since $\F$ is stable (this follows almost immediately
from the definition of stability and the fact that there are no non-trivial
division algebras over an algebraically closed field; see \cite[Cor.
1.2.8]{h-l}).

Similarly, an element in the kernel of $d_2^1$ comes from $H^1(C^{(p)},
\cEnd^0(\F))$, as asserted.
\end{proof}

What is interesting here is that not every
deformation $\E'$ of $\E$ arises in this way; indeed, $\E'$ admits a
$\nabla'$ if and only if it came from some deformation $\F'$ of $\F$. In particular, a given deformation of $\E$ admits a deformation of $\nabla$ if and only if it admits a deformation of $\nabla$ having $p$-curvature $0$.

We also get some information by computing the spectral sequence for the
double complex in the other direction. Specifically, we have:

\begin{prop}With the same notation and hypotheses as in Proposition
\ref{deg-les}, we have
$$H^1(C^{(p)}, \cEnd^0(\F)) \hookrightarrow \H^1(C, \K^\bullet)$$
The image is described by $1$-cocycles of $\cEnd^0(\E)$ in the kernel
of $d_\nabla$ (together with the zero $0$-cochain of $\cEnd^0(\E) \otimes
\Omega^1_C$).
\end{prop}

\begin{proof}Taking differentials in the opposite direction as before, we 
find our $E_2$ term looks like

$$\xymatrix{
{\vdots} & {\vdots} & {\vdots} & {}\\
{H^2(C, \ker d_\nabla)} & {H^2(C, \coker d_\nabla)} & {0} & {\dots}\\
{H^1(C, \ker d_\nabla)} & {H^1(C, \coker d_\nabla)} & {0} & {\dots}\\
{H^0(C, \ker d_\nabla)} & {H^0(C, \coker d_\nabla)} & {0} & {\dots}}$$

Now, $H^1(C, \ker d_\nabla)$ cannot have any further nonzero differentials
mapping into or out of it, so it injects into $\H^1(C, \K^\bullet)$.
We saw in the proof of the previous proposition that $\ker d_\nabla$
is precisely $F^{-1} \cEnd^0(\F)$. Since $F$ is a homeomorphism, we have 
$H^1(C, \ker d_\nabla) = H^1(C^{(p)}, \cEnd^0(\F))$, and we thus obtain the
desired injection, and description of its image. 
\end{proof}

The injectivity may be seen equally easily from the categorical equivalence between vector bundles $\F$ on $C^{(p)}$ and pairs $(\E, \nabla)$ on $C$.
Note, however, that this proposition does not imply that every non-trivial deformation of $\F$ maps to a non-trivial deformation of $\E$.

We now suppose that $\E$ is an extension of $\L^{-1}$ by $\L$, with $\L$ 
any line bundle on $C$. The short exact sequence
$$\xymatrix{
{0} \ar[r] & {\L} \ar[r]^-{i} & {\E} \ar[r]^-{j} & {\L^{-1}} \ar[r] & {0}
}$$
induces maps 
\begin{equation}d_1: \Hom(\L^{-1}, \L) \rightarrow \cEnd^0(\E)
\end{equation}
and
\begin{equation}\label{deg-d2} d_2: \cEnd^0(\E) \rightarrow \Hom(\L, \L^{-1})
\end{equation} 
by composition.
Specifically, $d_1(\phi)= i \circ \phi \circ j$, and $d_2(\phi) = j 
\circ \phi \circ i$. Thus, we obtain natural candidates for maps $H^0(C,
\Hom(\L^{-1}, \L)\otimes \Omega^1_C) \rightarrow \H^1(C, \K^\bullet)$ and
$\H^1(C, \K^\bullet) \rightarrow H^1(C, \Hom(\L,\L^{-1})$, although we will
need to check that the first is well-defined. We now restrict to our
particular situation of interest, for the remainder of this section and the
next.

\begin{sit}\label{deg-sit} Suppose that $C$ has genus $2$, and $\E$ an 
extension of $\L^{-1}$ by $\L$, with $\L$ a theta characteristic. Suppose 
also that the connection $\nabla$ on $\E$ has vanishing $p$-curvature.
\end{sit}

\begin{prop}In our situation, the maps $d_1$ and $d_2$ induce a short 
exact sequence
$$0 \rightarrow \Gamma(\cHom(\L^{-1},\L)\otimes \Omega^1_C) \rightarrow 
\H^1(C, \K^{\bullet}) \rightarrow \Ext(\L, \L^{-1}) \rightarrow 0$$
\end{prop}

\begin{proof} It is convenient to introduce the filtration of 
$\cEnd^0(\E)$ given by 
$0 \subset \im d_1 \subset \ker d_2 \subset \cEnd^0(\E)$.
We claim this induces the following filtration of our original complex:
\begin{equation}\label{deg-filt}
\xymatrix{
{\cEnd^0(\E)} \ar[r]^-{d_\nabla} & {\cEnd^0(\E)\otimes \Omega^1_C} \\
{\ker d_2} \ar@{^{(}->}[u] \ar[r]^-{d_\nabla} & {\cEnd^0(\E)\otimes \Omega^1_C} 
\ar@{^{(}->}[u] \\ 
{\im d_1} \ar@{^{(}->}[u] \ar[r]^-{d_\nabla} & {(\ker d_2)\otimes \Omega^1_C} 
\ar@{^{(}->}[u] \\
{0} \ar@{^{(}->}[u] \ar[r]^-{d_\nabla} & {(\im d_1) \otimes \Omega^1_C} \ar@{^{(}->}[u] 
}
\end{equation}

The only part which requires any verification is that 
$d_\nabla(\im d_1) \subset (\ker d_2) \otimes \Omega^1_C$, which is easily checked directly. Given this filtration of our complex, we get
(see \cite[1.4.5]{de1})
a spectral sequence converging to its hypercohomology, whose $E_2$ 
term is given in terms of the hypercohomology of the quotient complexes.
Specifically, if $\K^{\bullet}$ is our complex, and $\sA_i^{\bullet}$ the
filtration, we get $E_2^{p,q} = \H^{p+q}(C, \Gr^q_\sA (\K)) \Rightarrow
\H^{p+q}(C, \K)$. We must therefore start by calculating the associated
graded complexes of the filtration. 
Now, we have $\im d_1 \cong \cHom(\L^{-1}, \L)$, $\ker d_2/\im d_1 \cong \O_C$, and $\cEnd^0(E)/\ker d_2 \cong 
\cHom(\L, \L^{-1})$. Of course, the associated graded sheaves on the
right are gotten by taking these, shifted them down by one, and tensoring
with $\Omega^1_C$. Noting that $\cHom(\L^{-1}, \L) \cong \L^{\otimes 2} \cong 
\Omega^1_C$, we actually get isomorphic line bundles for the middle two
associated graded complexes, and we have to check that the maps between
them are isomorphisms. One can carry this out directly in terms of transition matrices, by checking that the two maps are linear and non-zero, at least in odd characteristic.

Noting that the hypercohomology of an isomorphism vanishes, and the
hypercohomology of a 1-term complex concentrated in the $i$th place is the
cohomology of the nonzero term, shifted by $i$, we obtain the
$E_2$ term of our spectral sequence:

$$\xymatrix@-10pt{
{H^0(C, \cH)} & {H^1(C, \cH)} & {0} & \dots & {} & {} \\
{} & {0} & {0} & {0} & {\dots} & {} \\
{} & {} & {0} & {0} & {0} & {\dots} \\
{} & {} & {} & {0} & {H^0(C, \cH\otimes \Omega^1_C)} & 
{H^1(C, \cH \otimes \Omega^1_C)}}$$
where $\cH := \cHom(\L, \L^{-1})$.

Now, there are potentially non-zero differentials between the remaining terms,
but because $H^0(C, \cHom(\L, \L^{-1}))=0$, the $H^0(C, \cH \otimes
\Omega^1_C)$ cannot have any further nonzero differentials, and applying $H^1(C, \cHom(\L, \L^{-1})) = \Ext(\L, \L^{-1})$, we obtain the
desired short exact sequence for $\H^1(C, \K^{\bullet})$.
\end{proof}

We also note from the construction of the spectral sequence that the
map from $\H^1(C, \K^{\bullet})$ to $\Ext(\L, \L^{-1})$ is, at least up
to sign, the map induced by first mapping to $H^1(C, \cEnd^0(\E))$, and 
then taking the map induced by $d_2$ on $H^1$. We thus get a diagram:

$$\xymatrix{
{} & {H^0(C, \cEnd^0(\E)\otimes \Omega^1_C)/H^0(C, \cEnd^0(\E))} 
\ar@{^{(}->}[d] \ar@{^{(}->}[dl] & {} \\
{H^0(C, (\Omega^1_C)^{2})} \ar@{^{(}->}[r] & {\H^1(C, \K^{\bullet})}
\ar[d] \ar@{->>}[r] & {\Ext(\L, \L^{-1})} \\
{} & {H^1(C, \cEnd^0(\E))} \ar[ur] & {}}$$
where the inclusion on the upper left follows formally from exactness in
the middle. But now we compare the dimensions to conclude:

\begin{prop}\label{deg-filtss}The kernel of $\H^1(C, \K^{\bullet}) \rightarrow H^1(C,
\cEnd^0(\E))$ is equal to the kernel of $\H^1(C, \K^{\bullet}) \rightarrow
\Ext(\L, \L^{-1})$. Equivalently, a deformation of the pair $(\E, \nabla)$
induces the trivial extension in $\Ext(\L, \L^{-1})$ if and only if it
was the trivial deformation of $\E$ (together with any deformation of
$\nabla$).
\end{prop}

\begin{proof} We need only show that the dimensions of the two kernels are
equal. Applying the Riemann-Roch theorem for vector bundles together with the self-duality of $\cEnd^0(\E)$, we compute that both 
$$h^0(C, (\Omega^1_C)^{\otimes 2})=h^0(C, \cEnd^0(\E)\otimes \Omega^1_C
)-h^0(C, \cEnd^0(\E))=3,$$
as desired.
\end{proof}

\section{Geometric Significance}\label{s-deg-geometric}

In this section, we lend some geometric substance to what we have calculated thus far, largely by characterizing the map $H^1(C, \cEnd^0(\E))
\rightarrow \Ext(\L, \L^{-1})$ in geometric terms. We continue with the
hypotheses of Situation \ref{deg-sit}. Suppose we have a family
of vector bundles $\tilde{\E}$ on $C$ with base $T$, and trivial determinant, 
such that for some $k$-valued point $0 \in T$, $\tilde{\E}|_0 \cong \E$. 
Since we have a map $\E \twoheadrightarrow \L^{-1}$, by adjointness we get a 
map $\tilde{\E} \twoheadrightarrow i_{0*} \L^{-1}$ (see Lemma
\ref{deg-adjoint}), and we can take the kernel to get 
a new family $\tilde{\E'}$ over $T$ which is isomorphic to $\tilde{\E}$ away 
from $0$ (this will be another vector bundle if $T$ is a reduced curve, but not
quite, for instance, if $T$ is $\Spec k[\epsilon]$; see Theorem \ref{deg-kerlf} 
of the appendix). 

We can then restrict back to the fiber at $0$, where we will get a new
$\E'$ on $C$ which will also be a vector bundle (even if $T = \Spec
k[\epsilon]$), of rank $2$ and trivial determinant. Everything but the 
trivial determinant assertion actually follows immediately from
Theorem \ref{deg-kerlf}, with $T$ either a curve or $\Spec k[\epsilon]$. For
the triviality of the determinant, and consequent remarks, we will for the
moment assume that
$T$ is a curve. In this case, we get the desired result from \cite[Cor
5.6]{mu1}, since we have triviality of the determinant away from $0$ on $T$.

Moreover, we see that $\E'$ is an extension of $\L$ by $\L^{-1}$:
Since $i$ is a closed immersion, $(i_{0*} \L^{-1})|_{0} = \L^{-1}$, so
we have an exact sequence 
\begin{equation}\label{deg-eel}\E' \rightarrow \E \rightarrow \L^{-1}
\rightarrow 0
\end{equation}
but since the kernel of $\E \rightarrow \L^{-1}$ is $\L$, $\E' \rightarrow
\E$ factors through $\L \rightarrow \E$, and we see that we get a
map from $\E'$ to $\L$, which is surjective, by exactness of equation 
\ref{deg-eel}. Since $\E'$ has trivial determinant, it follows that $\E'$
is an extension of $\L$ by $\L^{-1}$. Thus, in the case that $T$ is a curve,
we get a map

$$\phi_\E : \{\tilde{\E} \text{ over }T\} \rightarrow \Ext(\L, \L^{-1}),$$

This will have the following significance for us in trying to understand
the Verschiebung:

\begin{lem}\label{deg-efnot}Suppose $T$ is a curve, $0$ a point of $T$, and 
$\tilde{\F}$ a nontrivial family of 
semistable vector bundles over $T$ with trivial determinant, such that 
$F^* \F = \E$, where $\F := \tilde{\F}|_{0}$. That is to say, $\tilde{\F}$
gives a nonconstant map of $T$ into the moduli space $M_{2}$, passing
through a point where $V_{2}$ is undefined. Then writing $\tilde{\E} =
F^* \tilde{\F}$, if $\E' = \phi_{\E} (\tilde{\E}) \neq 0$, the limit point
of the image of $T$ under $V_{2}$ at $\F$ is given by $\E'$.
\end{lem}

\begin{proof} Since $\tilde{\E'}$ is isomorphic to $\tilde{\E}$ away
from the limit point, it suffices to check that $\E'$ is semistable as 
long as it is nonzero in $\Ext(\L, \L^{-1})$. But this may be checked directly from the definitions, using the fact that $\L$ has degree 1.
\end{proof}

The implication is that if we are lucky, we will be able to describe the
image of the exceptional divisor of $V_{2}$ in terms of $\P \Ext(\L,
\L^{-1})$, which can easily be checked to define a hyperplane inside of $M_{2}$, using the identification $M_2 \cong \P H^0(\Pic^1(C), \O(2 \Theta))$, together with \cite[Lem. 5.8]{n-r} to handle the semi-stable locus.
To show that this hyperplane will in fact be the image of the exceptional divisor, we first note that if $T = \Spec
k[\epsilon]$, we actually still get a map $\phi_\E : \Def^0(\E) \rightarrow
\Ext(\L, \L^{-1})$. In fact, more specifically, we have:

\begin{lem}\label{deg-samemap}$\phi_\E$ still exists in the case $T = \Spec
k[\epsilon]$, inducing a map from $\Def^0(\E)$ to $\Ext(\L, \L^{-1})$ 
which arises as the negative of the map on $1$-cocycles induced 
by $d_2: \cEnd^0(\E) \rightarrow \Hom(\L, \L^{-1})$ (\ref{deg-d2}), from the 
previous section.
\end{lem}

\begin{proof} Let $U_1, U_2$ be an open cover of $C$ trivializing $\L$, $\E$, and $\Omega^1_C$, and $\omega_i$ one-forms
trivializing $\Omega^1_C$ on the $U_i$. In addition, we set the convention that all 
$1$-cocycles will be written
with coordinates on $U_2$. 
To prove the lemma, we have to start by pinning down the
identification of $\Ext(\L, \L^{-1})$ with $H^1(C, \Hom(\L, \L^{-1}))$.
We will think of $\Ext(\L, \L^{-1})$ as being described by
transition matrices of the form 
$F = \begin{bmatrix}\vp_{12} & 0 \\ f & \vp_{12}^{-1} \end{bmatrix}$, 
where $f$ is any regular section on $U_1 \cap U_2$. We first claim 
that rather than taking $f$ itself as our $1$-cocyle of $\Hom(\L,\L^{-1})$,
we will have to take $-\vp_{12} ^{-1} f$. Indeed, this may be checked directly by following through the standard identification: starting from an 
extension $\F \in \Ext(\L, \L^{-1})$, the corresponding $1$-cocycle of
$\Hom(\L, \L^{-1})$ is obtained by considering the exact sequence 
$$0 \rightarrow \Hom(\L, \L^{-1}) \rightarrow \Hom(\L, \F) \rightarrow
\Hom(\L, \L) \rightarrow 0,$$
and looking at the image in $H^1(C, \Hom(\L, \L^{-1}))$ of the
identity in $\Hom(\L, \L)$ under the boundary map. 

Next, suppose we have a deformation $\E'$ of $\E$ given by the transition 
matrix 
$$E(I+ \epsilon E') =
\begin{bmatrix}\vp_{12} + \epsilon e_{11} & \vp_{12}^{-2} + \epsilon e_{12}\\
\epsilon e_{21} & \vp_{12}^{-1} + \epsilon e_{22}
\end{bmatrix}$$

Recalling that $d_2$, by definition, took the
lower left coordinate of a matrix, and noting that the lower left coordinate of $E'$ will be $\vp_{12} e_{21}$, we see from the above that 
it suffices to show that $\phi_\E (E')$ is described by the transition 
matrix 
$\begin{bmatrix}\vp_{12} & 0 \\ e_{21} & \vp_{12}^{-1}
\end{bmatrix}$. 
But we can calculate $\phi_\E (E')$ directly, in this case.
Denote by $s_i, t_i$ our trivializing basis on $U_i$, in terms of which 
$E$ and $E'$ are written. Also write $u_i$ for the trivialization
of $\L^{-1}$ on $U_i$. This means we have $s_2 = (\vp_{12} +
\epsilon e_{11}) s_1 + \epsilon e_{21} t_1$, 
$t_2 = (\vp_{12}^{-2} + \epsilon e_{12}) s_1 + (\vp_{12}^{-1}+\epsilon
e_{22}) t_1$, and $u_2 = \vp_{12}^{-1} u_1$. Now, the induced map from 
$\E'$ to $\L^{-1}$ sends $a s_i + b \epsilon s_i + c t_i + d \epsilon t_i$ 
simply to $c u_i$. This means that its kernel is generated by 
$s_i, \epsilon t_i$ on $U_i$. 
Using the above formulas for $s_2, t_2$ in terms of $s_1, t_1$, we find
that this kernel has transition matrix
$$\begin{bmatrix}\vp_{12} + \epsilon e_{11} & \epsilon \vp_{12}^{-2} \\
e_{21} & \vp_{12}^{-1} \end{bmatrix}.$$ 
Upon modding out by $\epsilon$,
this gives precisely the desired form for $\phi_\E (E')$. The theory 
developed in the appendix, and particulary Corollary \ref{deg-nr-trans} and 
the subsequent discussion, justifies this calculation, even though the 
transition matrix is not unique, and the kernel itself (prior to restriction) 
is not characterized by it.
\end{proof}

We next address some compatibility statements on curves and deformations.

\begin{lem}In the same situation as Lemma \ref{deg-efnot},
if we write $\bar{\F}$ for the 
induced first-order deformation of $\F$ gotten via some closed immersion 
$t_\epsilon: \Spec k[\epsilon] \hookrightarrow T$ deforming the point 
$0 \in T$, and $\bar{\E} := F^* \bar{\F}$, we have 
$\bar{\E} \cong t_\epsilon ^* \tilde{\E}$. Further, 
$\phi _\E \bar{\E} = \phi _\E \tilde{\E}$.
\end{lem}

\begin{proof}The first statement follows immediately from the fact that
$\tilde{\F}$ is a sheaf on $C^{(p)}\times T$, and all we are saying is
that pullback under $t_\epsilon$ commutes with pullback under $F$. Since
$t_\epsilon$ only acts on $T$, and $F$ only acts on $C^{(p)}$, they commute. 

The second half of the lemma is just an application of Theorem \ref{deg-switch}
in the appendix to our specific situation.
\end{proof}

Putting together this lemma with Lemma \ref{deg-efnot}, we see:

\begin{thm}\label{deg-geometric} In the same situation as the previous lemma, if we suppose
$\phi _\E \circ F^*$ is injective on first-order deformations of $\F$,
then the limit point at $0$ of the image of $T \rightarrow M_{2}$ 
under the Verschiebung is given as $\phi_\E \circ F^* (\bar{\F})$. In
particular, all such limit points are contained in $\P \Ext (\L, \L^{-1})
\subset M_{2}$. Further, the Verschiebung only needs to be blown up
once at $\F$.
\end{thm}

\begin{proof}Everything but the last assertion follows directly from
the two lemmas. To show that the Verschiebung only needs to be blown
up once at $\F$, it suffices to know that every smooth curve through $\F$ 
has the limit point of its image under $V_{2}$ determined by its tangent
at $\F$. It suffices to show that every such curve has a corresponding
family of vector bundles on it, since we could then apply our result that it
suffices to look at the first order deformation induced by the family.
However, because all of our $\F$'s are in the stable locus of our moduli
space, the obstruction to a universal bundle is given by a Brauer
class; see \cite[Cor. 4.3.5]{h-l}, \cite[Cor. 2.5, p. 55, Prop. 0.9, p. 
16]{m-f-k}, and \cite[Thm. 11.7]{gr3} to conclude that the GIT quotient is
an \'etale principal $\PGL_n$-bundle even in positive characteristic, and then
apply \cite[Prop. 3.3.2]{ca1}. Because $k$ is algebraically closed, by 
Tsen's theorem \cite[Rem. 1.14]{f-k}, \cite[Cor. 1.10]{gr2}
we obtain the required bundles on smooth curves.
\end{proof}

Finally, we can draw some conclusions which have immediate consequences for
our understanding of the Verschiebung:

\begin{thm}\label{deg-equivs}Given $\F$ such that $F^* \F \cong \E$, let $E_\F$ 
be the exceptional divisor above $\F$ after blowing up $M_{2}$ to make the
Verschiebung a morphism. Then of the following, 
a) and b) are equivalent, and either implies c) and d):
\begin{alist}
\itm The scheme of connections with vanishing $p$-curvature on $\E$
is reduced at the point corresponding to $\F$.
\itm The map $\Def(\F) \rightarrow \Def(\E)$ induced by $F^*$ is
injective.
\itm The image of $E_\F$ under $V_{2}$ in $M_{2}$ is precisely
$\P\Ext(\L, \L^{-1})$.
\itm $V_{2}$ only needs to be blown up once at $\F$
\end{alist}
\end{thm}

\begin{proof} a) is equivalent to there not being any non-trivial
deformations of $\nabla$ which hold $\E$ fixed, and have $p$-curvature
$0$. The equivalence of this with the condition that there are no nontrivial
deformations of $\F$ which pull back to the trivial deformation of $\E$ follows from the categorical equivalence provided by the Cartier isomorphism; see \cite[Thm. 5.1]{ka1}.

On the other hand, b) is equivalent to the map $\Def(\F) \rightarrow
\Ext(\L, \L^{-1})$ being injective, by Proposition \ref{deg-filtss}, which
implies it is an isomorphism, since both spaces have dimension $3$ over $k$. 
Noting that Lemma \ref{deg-samemap} tells us our geometric and cocycles versions 
of this map are really the same up to sign, by the preceding theorem this
implies c), as desired.

Lastly, the fact that d) follows from these conditions also follows from
the preceding theorem, since as we just noted, b) gives us that
$\Def(\F) \hookrightarrow \Ext(\L, \L^{-1})$.
\end{proof}

We may now easily put everything together to prove our main theorem:

\begin{proof}[Proof of Theorem \ref{deg-main}]
The assertion of (i) that each Frobenius-unstable bundle corresponds to an 
undefined point will follow from (iii), since the image of the exceptional
divisor of a blow-up centered at such a point is not just a single point.
The degree statement then follows from (ii) by Corollary
\ref{deg-versch}. But now (ii) and (iii) follow from
the implications b) implies c) and d) of Theorem \ref{deg-equivs}.
\end{proof}

We conclude with some further questions:

\begin{ques} Are statements c) and d) of Theorem \ref{deg-equivs} in fact
equivalent to a) and b)?
\end{ques}

\begin{ques} Is the scheme of Frobenius-unstable bundles of rank two and
trivial determinant (more precisely, the scheme of transport-equivalence 
classes of connections with trivial determinant and vanishing $p$-curvature
on the appropriate unstable bundles) isomorphic to the scheme-theoretic
undefined locus of $V_{2}$?
\end{ques}

\begin{ques} Is the degree of $V_{2}$ constant over all smooth curves of
genus 2?
\end{ques}

We remark that an affirmative answer to the second question would give an
affirmative answer to the third, thanks to a result of Mochizuki 
\cite[II, Thm. 2.8, p. 153]{mo3}: this gives that our scheme of 
Frobenius-unstable bundles is finite flat over our space of curves, and 
that it is smooth over the base field, from which it follows that its 
fiber over any fixed curve is a local complete intersection.

More generally, one could ask:

\begin{ques} How might one attempt to compute the degree of the Verschiebung
for curves of higher genus, or vector bundles of higher rank?
\end{ques}
 
We remark finally that to attempt to address this last question via 
similar techniques to those of this chapter, it would be necessary not only 
to generalize our understanding of the undefined locus of the Verschiebung, 
but also to appropriately generalize the degree formula of Proposition
\ref{deg-formula} to a substantially more general class of projective
varieties; indeed, for higher genus and rank, the moduli spaces in question
become singular along the strictly semi-stable locus.

\appendix

\section{Some General Results on the Verschiebung}\label{s-deg-background}

This appendix consists of the formal construction of, and some general
results on, the generalized Verschiebung map on coarse moduli spaces of 
vector bundles. Most of the arguments are due to A. J. de Jong, or Christian
Pauly, as indicated. We will work in the situation:

\begin{sit} $C$ is a smooth, proper curve over an algebraically closed field $k$ of characteristic
$p$. $C^{(p)}$ is the $p$-twist of $C$ over $k$, and $F$ is the relative 
Frobenius morphism from $C$ to $C^{(p)}$.
\end{sit}

Note that in this appendix, we do not assume $p>2$.

\begin{rem}The algebraically closed hypothesis is not actually necessary;
the moduli space construction can be made to work over a non-algebraically
closed field, and the arguments here will go through in this setting.
However, it is harder to find references for the general case, and we will
only apply the results here in the case of an algebraically closed base
field.
\end{rem}

\begin{rem} There are a few points to be careful of with respect to
characteristic and the general theory of moduli of vector bundles. The main 
obstruction is boundedness, which is not a problem in our situation of the
base being a curve (see \cite[Cor. 1.7.7]{h-l}), and is now known in any
dimension via the more involved argument of \cite{la2}. However, in
characteristic $p$ the statement that the moduli space universally
corepresents the relevant functor is in fact no longer true. It is however
true that it {\bf uniformly} corepresents the functor, which is to say that 
it is universal for flat base change, and this is all we will need; see
\cite[Thm. 1.10, p. 38]{m-f-k}.
\end{rem}

We recall the following definition and theorem:

\begin{defn}\label{deg-s-equiv} Two vector bundles $\E, \E'$ of degree $0$
are {\bf $S$-equivalent} if and only if there are filtrations $F_{\E}$
and $F_{\E'}$ with $\Gr_{F_{\E}} \cong \Gr_{F_{\E'}}$ (this isomorphism is
not required to preserve the grading), and the quotients of $F_{\E}$ and
$F_{\E'}$ all stable sheaves of degree $0$. 
\end{defn}

\begin{thm} \cite[Thm. 4.3.4]{h-l} There is a coarse moduli scheme $M_{r}(C)$ which uniformly
corepresents the functor of semistable vector bundles on $C$ of rank $r$ and 
trivial determinant. The closed points of $M_{r}(C)$ correspond to
$S$-equivalence classes of vector bundles; in particular, there is an open
subscheme $M^s_{r}(C)$ whose closed points parametrize stable vector
bundles.
\end{thm}

The main results on the Verschiebung are the following. The main arguments 
for existence and dominance are drawn from unpublished work of de Jong, 
while Pauly suggested the first portion of the argument for showing that the
map doesn't extend over the locus of Frobenius-unstable bundles.

\begin{thm}\label{deg-background} In the above situation, and given $r >0$, the operation of 
pulling back vector bundles under $F$ induces a {\bf generalized
Verschiebung} rational map $V_{r}: M_{r}(C^{(p)}) \dashrightarrow
M_{r}(C)$. If we denote by $U_{r}$ the open subset of
$M_{r}(C^{(p)})$ corresponding to bundles $\E$ such that $F^*(\E)$ is
semi-stable, we have further:
\begin{ilist}
\itm The domain of definition of $V_{r}$ is precisely $U_{r}$;
\itm $V_{r}$ is dominant.
\end{ilist}
\end{thm}

The existence portion of the theorem is straightforward.

\begin{prop} The $V_r$ of Theorem \ref{deg-background} exists, and induces a morphism on $U_r$.
\end{prop}

\begin{proof} Since semistability is an open condition (see \cite[Prop. 2.3.1]{h-l}), we get an open subfunctor $\U_{r}$ of the moduli 
functor $\cM_{r}(C^{(p)})$ corresponding to
semi-stable vector bundles on $C^{(p)}$ which pull back under $F$ to 
semi-stable vector bundles on $C$. We claim that it is enough to show that
this subfunctor is stable under $S$-equivalence. Indeed, given this, and using the description of the closed points of $M_{r}(C^{(p)})$ and the fact that it uniformly corepresents $\cM_{r}(C^{(p)})$, it is straightforward to check that $\U_{r}$ corresponds naturally to an
open subscheme $U_{r}$ of $M_{r}(C^{(p)})$ which corepresents $\U_{r}$ 
and whose closed points are precisely $S$-equivalence classes of vector 
bundles on $C^{(p)}$ whose pullbacks under $F$ are semi-stable. Now, 
Frobenius pullback induces a map 
from $\U_{r}$ to $\cM_{r}(C)$; if we compose with the map 
$\cM_{r}(C) \rightarrow M_{r}(C)$, since $U_{r}$ corepresents 
$\U_{r}$, we obtain the desired morphism $V_{r}: U_{r} \rightarrow 
M_{r}(C)$. 

We therefore show that our subfunctor is in fact stable under
$S$-equivalence. Let $\E, \E'$ be as in Definition \ref{deg-s-equiv}.
Since $C$ is smooth, we know $F$ is flat, so $F^*$ behaves well with
respect to filtrations and the operation $\Gr$, and we claim this implies
that $F^* \E$ is semi-stable if and only if $F^* \Gr_{F_{\E}}$ is
semi-stable. Certainly, if $\F \subset F^* \E$ is a destabilizing subsheaf,
then by considering the smallest subsheaf in the filtration $F^* F_{\E}$
into which $\F$ maps, there will be a nonzero map of $\F$ into the
corresponding quotient in $F^* \Gr_{F_{\E}}$. Conversely, suppose that $\F$
is a destabilizing subsheaf of some quotient in $F^* \Gr_{F_{\E}}$; that is,
$\F$ has positive degree and there is an injection into $F^* F_{\E}^{(i)}/F^*
F_{\E}^{(i+1)}$ for some $i$, from which it follows that the cokernel has
negative degree. If $\F'$ is the subsheaf of $F^* F_{\E}^{(i)}$
generated by $\F$ and $F^* F_{\E}^{(i+1)}$, the cokernel of the inclusion is
the same, hence has negative degree, and it follows that $\F'$ has positive
degree and maps into $F^* \E$, giving the desired instability.
We conclude that $F^* \E$ is semi-stable if and only if $F^*
\E'$ is semi-stable, as desired.
\end{proof}

Before completing the proof of the theorem, we recall a number of
well-known facts about moduli spaces of vector bundles.

\begin{thm} Let $\M_r(C)$ and $\M_r^{ss}(C)$ denote the Artin stacks of
vector bundles of rank $r$ and trivial determinant on $C$ and the open
substack of semistable bundles respectively. We have:
\begin{ilist}
\itm The Picard groups of $\M_r(C)$ and $\M_r^{ss}(C)$ are both isomorphic 
to $\Z$, and generated by the inverse of the determinant line bundle, which 
we denote by $\L$.
\itm For $n$ sufficiently large, $\L^n$ has base locus precisely equal to
the locus of unstable bundles, and descends to a very ample line bundle on 
$M_r(C)$.
\itm $M_r(C)$ is normal, and $\Q$-factorial.
\end{ilist}
\end{thm}

\begin{proof} To see that $\Pic(\M_r(C)) \cong \Pic(\M_r^{ss}(C))$ under the
natural map, use the argument of \cite[Prop. 8.3]{b-l}. Faltings gives a 
characteristic-independent argument for the isomorphism with $\Z$ in 
\cite[Thm. 17]{fa1}. 

The statement of (ii) in the context of GIT is essentially standard. 
Specifically, let $\bar{R}$ be a GIT rigidification space whose quotient 
is $M_r(C)$ \cite[\S 4.3]{h-l}, with linearized bundle $\hat{\L}$, and we 
will denote by $\hat{R}$ the open subscheme (containing $R$) of locally 
free sheaves which are globally generated with vanishing $H^1$, after 
twisting by the power of the ample bundle on $C$ used in the chosen 
rigidification. We see from vanishing of obstructions on curves and the
vanishing $H^1$ condition that $\hat{R}$ is smooth. The GIT-semistable 
points of $\hat{R}$ are precisely the semistable bundles 
\cite[Thm. 4.3.3]{h-l}, and it follows from the 
definition \cite[Def. 1.7 (b), p. 36, (2), p. 37]{m-f-k} that the rest are 
the points such that invariant sections of arbitrarily high powers of 
$\hat{\L}$ vanish; i.e., these form the base locus of high powers of 
$\hat{\L}$. Next, we claim that every unstable
bundle appears in some $\hat{R}$. But any bundle $\E$ generated by its 
global sections fits into an extension of $\det \E$ by $\O_C^{\oplus r-1}$
\cite[Thm. 2, p. 426]{at1}. Thus, for sufficiently high $n$, we find that
all $\E$ either semistable or equal to any chosen unstable $\E_0$ will be 
contained in the family
$$0 \rightarrow \O(-n)^{\oplus r-1} \rightarrow \E \rightarrow (\det \E)\otimes
\O(n(r-1))\rightarrow 0,$$ 
simply by ensuring that $n$ is high enough so that $\E(n)$ is globally
generated. Fixing $n$, if we choose $n'>>n$ high enough so that $\E(n')$
is globally generated and also $H^1(C, \E(n'))=0$ for all $\E$ in the above
$\Ext$ family, we find we can map this entire family into a corresponding
GIT $\Quot$ scheme, with the semistable locus mapping into $R$. But since 
the family is irreducible, it maps into $\bar{R}$, and we see from the 
definition of $\hat{R}$ that in fact it maps into $\hat{R}$, as desired. 
Finally, we need to argue that the stack statement follows. We remark that 
any $\hat{R}$ has quotient stack isomorphic to
an open substack of $\M_r(C)$ containing $\M_r^{ss}(C)$, and that 
furthermore as we vary our choices of $\hat{R}$, we obtain a cover of 
$\M_r(C)$. Also, using \cite[Lem. 8.2]{b-l} we see that $G$-invariant line 
bundles on $\hat{R}$ are in bijection with line bundles on $\M_r(C)$,
so $\hat{\L}$ must be a power of the pullback of $\L$. Now, the sections 
of $\L^{\otimes n}$ on $\M_r(C)$ are in bijection with the
sections on $\M^{ss}_r(C)$ \cite[Prop. 8.3]{b-l}, which are in bijection 
with $G$-invariant sections of the pullback on $\bar{R}^{ss}=R$ 
\cite[Lem. 7.2]{b-l}, 
which are in turn in bijection with the $G$-invariant sections of the
pullback on $\hat{R}$
\cite[Lem. 4.15]{n-r2}. Thus, every $G$-invariant section of the pullback
of $\L$ on $\hat{R}$ is
the pullback of a section on $\M_r(C)$, and since the varying choices of
$\hat{R}$ cover $\M_r(C)$, the base locus of $\L^{\otimes n}$ on $\M_r(C)$
is the union of the (images of the) base loci on the different $\hat{R}$,
which is precisely the unstable bundles, as desired.

For (iii), normality is evident because the GIT quotient uses rings of 
invariants of regular, hence integrally closed rings. Next, given (i),
the argument of \cite[9.2]{l-s1} shows that $M_r(C)$ is $\Q$-factorial.
\end{proof}

We now finish the proof of the theorem.

\begin{proof}[Proof of Theorem \ref{deg-background}] We begin by proving part (i), by showing that $V_{r}$ cannot be extended over the complement of $U_{r}$. The argument is in two parts. First, we show that if we know that the complement of $U_{r}$ has codimension at least $2$, then we obtain the desired result. 
We follow the notation of the preceding background theorem. Choose $n$ 
sufficiently large, so that $\L^{\otimes n}$ and
$(\L^{(p)}) ^{\otimes n}$ have base locus precisely equal to the locus of
unstable bundles, and descend to very ample line bundles $\bar{\L}$ and
$\bar{\L}^{(p)}$ on $M_r(C)$ and $M_r(C^{(p)})$
Let $(\bar{\L}, \bar{V})$ be the complete linear series for
$\bar{\L}$ on $M_r(C)$. 
We wish to show that $V_r^{*}(\bar{\L}, \bar{V})$
extends uniquely to a linear series $((\bar{\L}^{(p)})^{\otimes p},
V^{(p)})$ on $M_r(C^{(p)})$, whose base locus is precisely the locus of
Frobenius-unstable bundles. The uniqueness of this extension would then
imply that $V_r$ cannot be extended over any of the Frobenius-unstable
locus, as desired. Now, uniqueness follows from normality of $M_r(C^{(p)})$,
together with the codimension hypothesis on the complement of $U_r$. 
To see existence, we pull back to $\M_r^{ss}(C)$, where
by \cite[Prop. 8.3]{b-l}, we can extend to a linear series 
$(\L^{\otimes n}, V)$ on all
of $\M_r(C)$, which must have base locus precisely equal to the unstable locus. Using
\cite[Prop. A.1]{l-p2}, we pull this back under the Verschiebung to
$((\L^{(p)})^{\otimes np}, V^{(p)})$ on $\M_r(C^{(p)})$, which has base
locus precisely equal to the locus of bundles whose Frobenius-pullback are
unstable. Finally, by \cite[Prop. 8.4]{b-l} we can descend this to the
desired $((\bar{\L}^{(p)})^{\otimes p}, V^{(p)})$.

To complete the proof of (i), we show that the complement of $U_r$ has codimension at least $2$. Now, let $Z$ be any codimension $1$ closed subscheme of $M_r(C^{(p)})$; we claim that it necessarily intersects any positive-dimensional closed subscheme of $M_r(C^{(p)})$. Indeed, using the fact that $M_r(C^{(p)})$ 
is $\Q$-factorial, some multiple of $Z$ is necessarily Cartier, and in 
particular there is an effective Cartier divisor with set-theoretic support 
equal to $Z$. Moreover, because $\Pic(\M_r(C^{(p)}))\cong \Z$, this Cartier 
divisor is necessarily a descent to $M_r(C^{(p)})$ of some multiple of 
$\L^{(p)}$, which is ample, so must intersect any positive-dimensional 
closed subscheme, as claimed. But now we are done; the locus of direct 
sums of $r$ line bundles inside $M_r(C^{(p)})$ gives a $g^{r-1}$-dimensional closed subscheme of $M_r(C^{(p)})$ which is visibly contained in $U_{r}$, so we conclude that the codimension of $U_{r}$ is at least $2$, as desired.

We now prove part (ii). We simply exhibit a point of $M_{n}(C^{(p)})$ at
which $V_{n}$ induces a finite flat map on versal deformation spaces. This
point will correspond to a vector bundle $\E_0$ of the form $\L_1 \oplus 
\L_2 \oplus \cdots \oplus \L_n$, where the $\L_i$ are distinct line bundles 
of degree $0$, with $\bigotimes _i \L_i \cong \O_{C^{(p)}}$, so that 
$\E_0$ has trivial determinant. By \cite[p. 119]{ra2}, and in particular 
the asserted vanishing of $h^0(B \otimes L_1)$ for general $L_1$ and sequence 
(3), we may also require that the natural maps
$$H^1(C^{(p)}, \L_i \otimes \L_j^{-1}) \overset{F^*}{\rightarrow} H^1 (C,
F^* \L_i \otimes F^* \L_j)$$
are isomorphisms for all $i \neq j$.

We will consider deformations of $\E_0$ preserving the triviality of the
determinant. One can check that this satisfies the criteria $(\oH_1)$, 
$(\oH_2)$, $(\oH_3)$ of \cite{sc2}, so we let $R^{(p)}$ be a hull for 
the deformation problem and $\E$ on $C^{(p)} \times \Spec R^{(p)}$ be the
corresponding versal deformation of $\E_0$. The first-order deformations are parametrized by $H^1(C^{(p)}, \cEnd^0(\E_0))$, and deformations are 
visibly unobstructed, so we conclude that 
$R^{(p)} \cong k[[t_1, \dots t_N]]$, with 
$N = h^1(C^{(p)}, \cEnd^0(\E_0))$. Since 
$$\cEnd^0(\E_0) \cong \O_{C^{(p)}} ^{\oplus n-1} \oplus 
\bigoplus _{i \neq j} \L_i \otimes \L_j^{-1},$$
Riemann-Roch for vector bundles gives 
\begin{multline*}h^1(C^{(p)}, \cEnd^0(\E_0)) = h^0(C^{(p)}, \cEnd^0(\E_0)) 
- \deg \cEnd^0(\E_0) + (\rk \cEnd^0(\E_0))(g-1) \\
= (n-1)- 0 + (n^2-1)(g-1) = (g-1) n^2 + n - g.
\end{multline*}  

Now, let $I^{(p)} \subset R^{(p)}$ be the ideal defining the maximal closed
subscheme of $\Spec R^{(p)}$ over which $\E$ remains isomorphic to a direct 
sum of $n$ distinct line bundles. If $J^{(p)}$ denotes the Jacobian of 
$C^{(p)}$, and $(J^{(p)})^n \rightarrow J^{(p)}$ is the addition morphism, 
let $T^{(p)}$ be the fiber of this morphism over $0$. Since the kernel of 
the addition map corresponds precisely to $n$-tuples of line bundles whose 
direct sum has trivial determinant, we may 
then describe $R^{(p)}/I^{(p)}$ as the completion of the local ring of 
$T^{(p)}$ at the point corresponding to $\L_1, \dots, \L_n$. Now, if we set 
$$\F_0 := F^* \E_0 \cong \bigoplus _i F^* \L_i,$$
we can as before let $R$ be a hull for the deformations of $\F_0$, and $I$
the ideal cutting out the locus preserving the direct sum decomposition. We
then obtain the corresponding description of $R/I$ as above. Now, it is easily checked that the pullback under Frobenius induces a morphism 
$v: R \rightarrow R^{(p)}$. Next, since direct sum decompositions are 
preserved under pullback, $v(I) \subset I^{(p)}$, so we get an induced 
homomorphism $R/I \rightarrow R^{(p)}/I^{(p)}$. This clearly
corresponds to the morphism $T^{(p)} \rightarrow T$ induced by the
Verschiebung morphism $V: J^{(p)} \rightarrow J$; this last is finite flat,
so we find that $R/I \rightarrow R^{(p)}/I^{(p)}$ is also finite flat.

Now, we also assert that we have 
$$\Hom_k(I^{(p)}/ \m_{R^{(p)}} I^{(p)}, k) \cong 
\bigoplus_{i \neq j} H^1 (C^{(p)}, \L_i \otimes \L_j ^{-1})$$
and similarly for $\Hom_k(I/ \m_R I,k)$; we illustrate the argument for 
$R$. By the definition of a hull, $\m_R/\m_R^2$ is dual to the space
of first-order infinitesmal deformations $\T_{\F_0}$ of $\F_0$, with 
the subspace $\T'_{\F_0} \subset \T_{\F_0}$ obtained by modding out by
$I$ and then dualizing corresponding to deformations preserving the direct 
sum decomposition. Thus, the quotient $\T_{\F_0}/\T'_{\F_0}$ corresponds to 
deformations transverse to those preserving the direct sum decomposition, 
and is obtained as the dual of the subspace 
$I/\m_R I \subset \m_R / \m_R^2$. Now, it is also easy to see that 
$\T'_{\F_0}$ also corresponds to the summand of $H^1(C, \cEnd^0(\F_0))$ 
given by $H^1(C, \O_C^{\oplus n-1})$. Thus, we see that 
$\T_{\F_0}/\T'_{\F_0}$ is given by
$H^1 (C, \bigoplus_{i \neq j} F^* \L_i \otimes F^* \L_j^{-1}) 
= \bigoplus _{i \neq j} H^1(C, F^* \L_i \otimes \L_j^{-1})$, 
which is the desired statement for $R$.

As a result, we find by our additional hypotheses on the choice of the
$\L_i$ that $v$ induces an isomorphism 
$I/\m_R I \risom I^{(p)}/\m_{R^{(p)}} I^{(p)}$; we claim that this
together with the finite-flatness of $R/I \rightarrow R^{(p)}/I^{(p)}$ 
is enough to imply that $v: R \rightarrow R^{(p)}$ is itself 
flat. Indeed, the isomorphism $I/\m_R I \risom I^{(p)}/\m_{R^{(p)}} 
I^{(p)}$ implies that a set of generators of $I/\m_R I$ over
$R$ will map to a set of generators of $I^{(p)}/\m_{R^{(p)}} I^{(p)}$ over 
$R$ and in particular over $R^{(p)}$, which by Nakayama's lemma implies that 
any lifts generate $I^{(p)}$ over $R^{(p)}$; we conclude that $I^{(p)} =
v(I) R^{(p)}$. Now, we note that because $I_r \subset \m_R$, we may write
$R^{(p)}/v(\m_R) R^{(p)} = (R^{(p)}/v(I) R^{(p)})/v(\m_R)(R^{(p)}/v(I)
R^{(p)}) = (R^{(p)}/I^{(p)})/v(\m_R)(R^{(p)}/I^{(p)})$, which must be finite
over $R/I$ and hence over $R$ because $R^{(p)}/I^{(p)}$ is. We thus have
that the closed fiber of $\Spec R^{(p)} \rightarrow \Spec R$ is finite, and
since both rings are regular of the same dimension, by \cite[Thm 18.16
b]{ei1} the map is flat, hence dominant.
It remains to check that the map on hulls being dominant implies that 
the map on coarse moduli spaces is dominant. Because the coarse moduli
spaces are irreducible (see \cite[Rem. 5.5]{ne1}), it suffices to show 
that the image of the map from the hull to the coarse 
moduli space is not contained in any proper closed subset, which follows easily from the definitions.
\end{proof}

We conclude with a lemma on the case $g=2, r=2$ case which we used to 
conclude the Verschiebung is given by polynomials of degree $p$ in 
Section \ref{s-deg-prelim}.

\begin{lem}\label{det-o1} When $g=r=2$, the line bundle $\O(1)$ on 
$M_2(C)\cong \P^3$
pulls back to the inverse of the determinant bundle on the stack $\M_2(C)$.
\end{lem}

\begin{proof}We first use the criterion of \cite[Prop. 2.2 (B)]{m-r} to 
descend the determinant bundle $\L$ away from $\O_C \oplus \O_C$. On this 
locus, one easily check that the two conditions are satisfied for points 
with closed orbits, as these are the polystable bundles. Because $\P^3$ is 
regular, we can then extend uniquely to a descent of $\L$ on all of 
$M_2(C)$. Finally, since we already noted above that the inverse of the
determinant bundle is the ample generator of $\Pic(\M_r(C))$, the only 
bundle it could descend to is $\O(1)$. 
\end{proof}

\section{A Commutative Algebra Digression}\label{s-deg-digress}

We develop some
simple but non-standard commutative algebra over non-reduced rings which
will be helpful in studying deformations of vector bundles. Throughout
this section, $R$ will denote a (typically non-reduced) Noetherian ring,
$\n$ the ideal of nilpotents of $R$, and $M$ a finitely generated $R$-module, although we drop the Noetherian and finite-generated hypotheses for the first lemma.
Note that everything we define in this section will be equivalent to their
standard versions whenever $R$ is integral.

\begin{defn}$M$ is {\bf $NR$-free} of rank $r$ if $M$ is generated by some
$m_1, \dots m_r$, such that given any relation $\sum_i a_i m_i =0$, all
the $a_i$ must be nilpotent.
\end{defn}

\begin{lem}If $M$ is $NR$-free of rank $r$, then $M_{\red}$ is 
free of rank $r$. In particular, rank is well-defined for $NR$-free
modules. Further, the converse holds if $M$ is finitely generated and $R$ is local, or if $\n ^m =0$ for some $m$.
\end{lem}

\begin{proof}Suppose $M$ is $NR$-free of rank $r$. Clearly, the $m_i$
still generate $M_{\red}$, but using the fact that the $m_i$ generate $M$, one checks that any non-zero relation modulo the nilpotents of $R$ would yield a relation in $R$ with not all coefficients nilpotent. Conversely, suppose 
$M_{\red}$ is free of rank $r$, with generators $m_i$. We 
show that any lifts of the $m_i$ generate $M$, and they then clearly
satisfy the desired relations restrictions. If $R$ is local with
maximal ideal $\m$, and $M$ finitely generated, this follows immediately from Nakayama's lemma. If, on the other hand, $\n ^m =0$, we will get the
desired result by showing via induction that $M' := M/(\{m_i\}_i)$
is contained in (the image of) $\n^j M$ for all $j$. The base case is $j=0$, 
which is a triviality. Now suppose it holds for $j-1$, and we want to show
it for $j$. Take $m \in M$; by the induction hypothesis, we can write
$m = m' + \sum _i a_i m_i$ for some $m' \in \n^{j-1} M$. Since the $m_i$ generate $M_{\red}$, writing out $m'$ explicitly as a sum of products of elements of $\n$, we find that there exist $a'_i \in \n^{j-1}$ such that $m'-\sum_i a'_i m_i \in \n^j M$, giving the desired result.
\end{proof}

This argument also immediately gives us:

\begin{cor}\label{deg-nr-trans} Given two generating sets, $m_i$ and $m'_i$
for an $NR$-free
$R$-module $M$, there is a (non-unique) invertible matrix $T$
relating the $m_i$ to the $m'_i$, such that if $\bar{m}_i$, $\bar{m}'_i$,
and $\bar{T}$ are the images of $m_i$, $m'_i$ and $T$ in $M_{\red}$, 
$\bar{T}$ is the matrix relating the $\bar{m}_i$ to the
$\bar{m}'_i$.
\end{cor}

We can extend the standard definitions.

\begin{defn}$M$ is {\bf locally $NR$-free} of rank $r$ if $M$ becomes 
$NR$-free of rank $r$ over every local ring of $R$. Given a separated 
Noetherian scheme $X$, a coherent sheaf $\F$ of 
$\O_X$-modules is locally $NR$-free of rank $r$ if $\F(U)$ is locally
$NR$-free of rank $r$ over $\O_X(U)$ for every affine $U$. Because our 
modules
are finitely generated, this is equivalent to there being an open cover
of $X$ on which $\F$ becomes $NR$-free of rank $r$.
\end{defn}

We see that while we can attach transition matrices to locally $NR$-free
sheaves, they are not unique, nor do they uniquely determine $\F$. However,
they do determine $\F_{\red}$ over $X_{\red}$, which is all we will need for
our purposes.

To develop the results we want, we will make the hypothesis for the rest
of this section that $\Spec R$ is irreducible. With this, we can prove the
following proposition, which properly reformulated will kill several
birds with one stone:

\begin{prop}\label{deg-kernr}Suppose that $R_{\red}$ is integral, and let $M$ be a 
free $R$-module of rank $r$, and $N$ an $R$-module generated by $n_1, \dots, 
n_s$, with $\ell$ non-zero relations of the form $f n_i =0$ for some $f \in R$.
Then if $\phi$ is a surjective map from $M$ to $N$, $\ker \phi$ is 
an $NR$-free $R$-module of rank $r-s+\ell$. If further $\ell=0$, then $\ker \phi$
is in fact free.
\end{prop}

\begin{proof}First, we can suppose that $f$ is not a unit, as if it were
$N$ would just be a free module of rank $s-\ell$, and we could simply use
$n_{\ell+1}, \dots, n_s$ as generators. Our first claim is that if we lift
$\phi$ from $N$ to $\tilde{N} := R[n_1, \dots, n_s]$, the map will remain
surjective. Let $\tilde{\phi}$ be some such lift, one can check by induction on $j$ that $\tilde{N}
/\im \tilde{\phi} \subset f^j \tilde{N}/\im \tilde{\phi}$ for all $j$.
The Krull Intersection Theorem then implies that for some $f' \in (f)$, we 
have $(1-f')(\cap _j (f^j)) = 0$. Now, by the
hypothesis that $R_{\red}$ is integral, either $1-f'$ is nilpotent, or
$\cap _j (f^j) = 0$. However, if $1-f'$ were nilpotent, $f'$ would have to
be a unit, which is not possible, since $f$ was assumed not to be a unit.
It follows that $\cap _j (f^j)=0$, and $\tilde{\phi}$ is surjective.

Next, since it is a surjective map between free modules, the kernel of 
$\tilde{\phi}$ must be free, of rank $r-s$. This immediately proves the 
last assertion of the theorem. Now, let $\ker \tilde{\phi}$ be generated by
$\tilde{m}_1, \dots, \tilde{m}_{r-s}$. Next, let $\hat{m}_1, \dots, 
\hat{m}_{\ell}$ be any elements mapping to $f n_1, \dots, f n_{\ell} \in
\tilde{N}$ under $\tilde{\phi}$. 
It is then a routine check from the definitions $\{\tilde{m}_i, \hat{m}_j\}$ is a generating set making $\ker \phi$ into an $NR$-free module, of rank $r-s+\ell$.
\end{proof}

\begin{rem}Note in particular that if 
$f$ is a nonzero nilpotent, we will obtain modules of different rank
by first taking the kernel and then restricting to $R_{\red}$
than by restricting to $R_{\red}$ and then taking the kernel.
\end{rem}

To apply this proposition, we define:

\begin{defn}An {\bf effective $NR$-Cartier divisor} on a Noetherian, 
separated, irreducible scheme $X$ is a global section of the monoid sheaf
$(\O_X\smallsetminus\{0\})/\O_X^*$.
\end{defn}

\begin{lem}Associated to any non-trivial effective $NR$-Cartier divisor $f$ 
on $X$ is a canonical closed immersion $X_f \hookrightarrow X$. Given any 
closed subscheme of $X$, there is at most one effective $NR$-Cartier divisor
which induces it.
\end{lem}

\begin{proof}Clear.
\end{proof}

\begin{rem}Although $X_f$ need not have codimension $1$ in $X$ for an 
effective $NR$-Cartier, it behaves in certain ways as if it had codimension 
$1$, as the following theorem demonstrates.
\end{rem}

\begin{thm}\label{deg-kerlf}
Let $X_f \hookrightarrow X$ be the closed immersion associated
to a nontrivial effective $NR$-Cartier divisor $f$ on $X$, $\F$
a locally free $O_X$-module of rank $r$ on $X$, and $\sG$ a coherent 
$\O_X$-module on $X$. Let $f: \F \rightarrow \sG$ be surjective. Then:
\begin{ilist}
\itm if $\sG$ is locally free of rank $s$ on $X$, $\ker f$ is locally 
free of rank $r-s$, and
\itm if $\sG$ is the pushforward of a locally free sheaf of rank $s$
on $X_f$, $\ker f$ is locally $NR$-free on $X$ of rank $r$.
\end{ilist}
\end{thm}

\begin{proof}(i) follows immediately from Proposition \ref{deg-kernr},
by restricting to a cover on which $\F$ and $\sG$ are free, and noting
that we are in the case $\ell=0$. (ii) also follows from the same proposition,
as if our cover is also fine enough give trivializing elements for $f$,
we find that $f$ matches up with the $f$ of the proposition, and because
$\sG$ is a pushforward from $X_f$, it satisfies the hypotheses of the
proposition with $\ell=s$, giving us the desired conclusion.
\end{proof}

We will immediately apply this to prove a theorem to the effect that
in certain cases, when one wants to replace an unstable vector bundle in a
family with a semistable one, it suffices to look simply at the first order
deformation induced by the family. We begin with the following lemma.

\begin{lem}\label{deg-adjoint} Let $f: S \hookrightarrow T$ be a closed immersion of schemes,
$\F$ a locally free sheaf on $T$, and $\sG$ a locally free sheaf on
$S$. Given $\psi: \F \rightarrow f_* \sG$, then $\psi$ arises via adjointness
from $f^* \psi : f^* \F \rightarrow f^* f_* \sG = \sG$, and if $f^* \psi$
is surjective, then so is $\psi$. In this case, there is a natural surjection 
$f^* \ker \psi \twoheadrightarrow \ker f^* \psi$.
\end{lem}

\begin{proof} Routine.
\end{proof}

\begin{thm}\label{deg-switch} Let $f: S \hookrightarrow T$ be a closed 
immersion, 
with $T$ reduced, and such that if $i: S_{\red} \hookrightarrow S$ is the 
standard inclusion, then both $i$ and $f \circ i$ are induced by 
effective $NR$-Cartier divisors. Let $\F$ be a locally free
sheaf of rank $r$ on $T$, and $\sG$ a locally free sheaf of rank $s$ on
$S_{\red}$, and $\psi: \F \rightarrow f_* i_* \sG$ arising via adjointness 
from a surjective map $f^* \F \rightarrow i_* \sG$, which must then be 
$f^* \psi$. Then $i ^* f^* \ker \psi = i^* \ker f^* \psi$, and both are 
locally free of rank $r$ on $S_{\red}$.
\end{thm}

\begin{proof}The previous lemma gives us that $\psi$ is surjective, and
a surjection $f^* \ker \psi
\twoheadrightarrow \ker f^* \psi$. Since pullback is right exact,
this gives us $i^* f^* \ker \psi \twoheadrightarrow i^* \ker f^* \psi$,
and by Theorem \ref{deg-kerlf} (i), it suffices to show these are both
locally free of rank $r$, since then we would have the kernel of this
surjection also being locally free, of rank $0$. On the other hand, by (ii) 
of the same theorem, $\ker \psi$ is is locally $NR$-free on $T$ of rank 
$r$, hence locally free of rank $r$ by reducedness of $T$. Thus, 
$i^* f^* \ker \psi$ must also be locally
free of rank $r$. But likewise, since $f^* \F$ must also be locally
free of rank $r$, using the same theorem, $\ker f^* \psi$ is locally 
$NR$-free of rank $r$ on $S$, which as we saw initially is equivalent
to $i^* \ker f^* \psi$ being locally free of rank $r$ on $S_{\red}$,
as desired.
\end{proof}

\bibliographystyle{hamsplain}
\bibliography{hgen}
\end{document}